\documentclass[11pt]{article}
\usepackage[english]{babel}
\usepackage{amsmath, amsthm, fancyhdr}

\newcommand{\version}{version 1.0,\ \ Oct. 31, 2020}

\usepackage{amsmath, amsthm, fancyhdr}
\usepackage{amscd}
\usepackage{amsfonts}
\usepackage{amssymb}
\usepackage{mathrsfs}
\usepackage{url}
\usepackage[toc,page]{appendix}

\newcommand{\R}{\mathbb{R}}
\newcommand{\C}{\mathbb{C}}

\newcommand{\Z}{\mathbb{Z}}
\DeclareMathOperator{\NS}{NS}
\DeclareMathOperator{\FS}{FS}
\newcommand{\sO}{\mathscr{O}}

\DeclareMathOperator{\Gr}{Gr}
\DeclareMathOperator{\Deg}{Deg}

\DeclareMathOperator{\End}{End}
\DeclareMathOperator{\Id}{Id}
\newcommand{\1}{\sqrt{-1}}

\renewcommand{\bar}{\overline}
\renewcommand{\phi}{\varphi}
\renewcommand{\epsilon}{\varepsilon}
\renewcommand{\geq}{\geqslant}

\def\eqref#1{(\ref{#1})}

\newcommand{\rk}{{\rm rk}}
\newcommand{\arrow}{{\:\longrightarrow\:}}
\newcommand{\Aut}{\operatorname{Aut}}
\newcommand{\restrict}[1]{{\left|_{{\phantom{|}\!\!}_{#1}}\right.}}
\newcommand{\cntrct}                % contraction with a vector field
{\hspace{2pt}\raisebox{1pt}{\text{$\lrcorner$}}\hspace{2pt}}

\numberwithin{equation}{section}

\newcommand{\arr}{\xrightarrow}

\renewcommand{\o}{\otimes}

\theoremstyle{definition}

\newcounter{Mycounter}[section]
\newcounter{lemma}[section]
\renewcommand{\thelemma}{{Lemma \thesection.\arabic{lemma}}}
\newcommand{\lemma}{%
     \setcounter{lemma}{\value{Mycounter}}
     \refstepcounter{lemma}
     \stepcounter{Mycounter}
     {\bf \thelemma:\ }}

\newcounter{claim}[section]
\renewcommand{\theclaim}{{Claim \thesection.\arabic{claim}}}
\newcommand{\claim}{%
     \setcounter{claim}{\value{Mycounter}}
     \refstepcounter{claim}
     \stepcounter{Mycounter}
     {\bf \theclaim:\ }}

\newcounter{sublemma}[section]
\newcounter{corollary}[section]
\newcounter{theorem}[section]
\renewcommand{\thetheorem}{{Theorem \thesection.\arabic{theorem}}}
\newcommand{\theorem}{%
     \setcounter{theorem}{\value{Mycounter}}
     \refstepcounter{theorem}
     \stepcounter{Mycounter}
     {\bf \thetheorem:\ }}

\newcounter{conjecture}[section]
\newcounter{proposition}[section]
\renewcommand{\theproposition}
       {{Proposition \thesection.\arabic{proposition}}}
\newcommand{\proposition}{%
     \setcounter{proposition}{\value{Mycounter}}
     \refstepcounter{proposition}
     \stepcounter{Mycounter}
     {\bf \theproposition:\ }}

\newcounter{definition}[section]
\renewcommand{\thedefinition}
       {{Definition~\thesection.\arabic{definition}}}
\newcommand{\definition}{%
     \setcounter{definition}{\value{Mycounter}}
     \refstepcounter{definition}
     \stepcounter{Mycounter}
     {\bf \thedefinition:\ }}

\newcounter{example}[section]
\newcounter{remark}[section]
\renewcommand{\theremark}{{Remark \thesection.\arabic{remark}}}
\newcommand{\remark}{%
     \setcounter{remark}{\value{Mycounter}}
     \refstepcounter{remark}
     \stepcounter{Mycounter}
     {\bf \theremark:\ }}

\newcounter{problem}[section]
\newcounter{question}[section]
\makeatletter
\@addtoreset{equation}{section}
\@addtoreset{footnote}{section}
\makeatother

\def\blacksquare{\hbox{\vrule width 5pt height 5pt depth 0pt}}
\def\endproof{\blacksquare}
\def\proof{{{\bf Proof: }}}

\begin{document}

%%%%%%%%%%%%%%%%%%%%%%%%%%%%%%%%%%%%%%%%%%%%%%%%%%%%%%%%%%%%
\begin{center}
{\Large\bf
Sections of Lagrangian fibrations \\[2mm] on 
holomorphically  symplectic manifolds\\[2mm]
and degenerate twistorial deformations\\[4mm]}
Fedor A. Bogomolov\footnote{Fedor Bogomolov is 
partially supported by the HSE University Basic
Research Program, Russian Academic 
Excellence Project '5-100' and by EPSRC
programme  grant  EP/M024830.}, Rodion N. D\'eev, Misha Verbitsky\footnote{Partially supported by the  Russian Academic Excellence Project '5-100',
FAPERJ E-26/202.912/2018 and CNPq - Process 313608/2017-2.}

\end{center}
%%%%%%%%%%%%%%%%%%%%%%%%%%%%%%%%%%%%%%%%%%%%%%%%%%%%%%%%%%%%

%%%%%%%%%%%%%%%%%%%%%%%%%%%%%%%%%%%%%%%%%%%%%%%%
{\small 
\hspace{0.15\linewidth}
\begin{minipage}[t]{0.7\linewidth}
{\bf Abstract} \\
Let $(M,I, \Omega)$ be a holomorphically symplectic manifold
equipped with a holomorphic Lagrangian fibration $\pi:\; M \mapsto X$, and
$\eta$ a closed form of Hodge type (1,1)+(2,0) on $X$.
We prove that $\Omega':=\Omega+\pi^* \eta$ is again 
a holomorphically symplectic form, for another
complex structure $I'$, which is uniquely determined by $\Omega'$.
The corresponding deformation of complex structures
is called ``degenerate twistorial deformation''.
The map $\pi$ is holomorphic with respect to this
new complex structure, and $X$ and the fibers
of $\pi$ retain the same complex structure as
before. Let $s$ be a smooth section of of $\pi$. We prove that
there exists a degenerate twistorial deformation 
$(M,I', \Omega')$  such that $s$ is a holomorphic
section.
\end{minipage}
}

\tableofcontents

%%%%%%%%%%%%%%%%%%%%%%%%%%%%%%%%%%%%%%%%%%%%%%%%%%%%%%%%%%%%%%%%%%%%%

\section{Introduction}

%%%%%%%%%%%%%%%%%%%%%%%%%%%%%%%%%%%%%%%%%%%%%%%%%%%%%%%%%%%%%%%%%%%%%

%%%%%%%%%%%%%%%%%%%%%%%%%%%%%%%%%%%%%%%%%%%%%%%%%%%%%%%%%%%%
\subsection{Complex structure obtained from a 
complex symplectic form}
%%%%%%%%%%%%%%%%%%%%%%%%%%%%%%%%%%%%%%%%%%%%%%%%%%%%%%%%%%%%

A complex structure $I$ on a real vector space $V$
is uniquely determined by a complex-linear symplectic
form $\Omega$. Indeed, $\Omega$ has Hodge type $(2,0)$
because it is complex-linear. The corresponding
map $\Omega:\; V\otimes_\R \C \arrow V^*\otimes_\R \C$
is non-degenerate on the Hodge component $V^{1,0}$,
and vanishes on $V^{0,1}$. Since the complex structure
on $V$ is uniquely determined by $V^{0,1}\subset V\otimes_\R \C$,
$\Omega$ determines $I$ uniquely.

This description is very beneficial on a manifold,
when $\Omega$ is a complex-valued differential form.
In that case, $d\Omega=0$ implies that the almost
complex structure determined by $I$ is integrable
(\ref{clos}). In \cite{_Verbitsky:Degenerate_},
this observation was used to construct deformations
of hyperk\"ahler manifolds admitting a Lagrangian
fibration. It turns out that a holomorphic
symplectic form can be characterized intrinsically
in terms of its rank and exterior powers 
(\ref{_inhe_complex_stru_Proposition_}).

More precisely, a complex-valued exterior
2-form $\Omega$ on a $4n$-dimensional real vector space
is called {\bf c-symplectic} if $\Omega^{n+1}=0$
and $\Omega^n \wedge \bar\Omega^n$ is non-degenerate.
By \ref{_inhe_complex_stru_Proposition_}, this is
equivalent to being a non-degenerate complex linear
2-form for some complex structure on $V$.

Recall that {\bf holomorphic Lagrangian fibration}
on a holomorphically symplectic manifold $(M, \Omega)$
is a holomorphic map with fibers which are Lagrangian
with respect to $\Omega$.

Given a Lagrangian fibration $\pi:\; M \arrow B$ on a holomorphic
symplectic manifold $(M, \Omega)$ and a closed 2-form
$\eta$ on $B$ of Hodge type $(2,0)+(1,1)$, 
the sum $\Omega+\pi^*\eta$ is c-symplectic,
hence determines an almost complex structure
$I_\eta$ on $M$ (\ref{_inhe_complex_stru_Proposition_}).
This almost complex structure is integrable
when $\eta$ is closed (\ref{clos}).
The corresponding deformation of complex structures
is uniquely determined by the cohomology class
of $\eta$; it is called {\bf degenerate twistor
  deformation} (\cite{_Verbitsky:Degenerate_}).

As shown in \cite{_Verbitsky:Degenerate_},
the map $\pi:\; M \arrow B$ is a  holomorphic
Lagrangian fibration with respect to $\eta$,
and the complex structure on the fibers of $\pi$ 
does not change; in dimension 1, this deformation  
gives ``the Tate-Shafarevich twist'' of an elliptic
fibration (that is, a deformation of a fibration
which retains the complex structures on the base
and on the fibers). 

%%%%%%%%%%%%%%%%%%%%%%%%%%%%%%%%%%%%%%%%%%%%%%%%%%%%%%%%%%%%
\subsection{Lagrangian fibrations on hyperk\"ahler
  manifolds}
\label{_Lagra_fibra_Subsection_}
%%%%%%%%%%%%%%%%%%%%%%%%%%%%%%%%%%%%%%%%%%%%%%%%%%%%%%%%%%%%

When $M$ is a compact hyperk\"ahler manifold
of maximal holonomy,\footnote{This is the same
as ``IHS'', irreducibly holomorphic symplectic manifold.}
 essentially all non-trivial
holomorphic maps $M\arrow B$ are Lagrangian fibrations,
by Matsushita's theorem. 

\hfill

\theorem
Let $M$ be a compact hyperk\"ahler manifold
of maximal holonomy, and $\pi:\; M\arrow B$
a surjective holomorphic map, with $\dim M > \dim B >0$.
Then $\pi$ is a holomorphic Lagrangian fibration.

\proof \cite{_Matsushita:fibred_}. \endproof

\hfill

Applying the Stein factorization theorem, we
can factorize a given holomorphic Lagrangian fibration
through a fibration with connected fibers. 
Further on, we shall tacitly assume that all
Lagrangian fibrations we consider have 
connected fibers. In this case the base 
$B$ is known to have
the same rational cohomology as a complex projective space,
by a theorem of D. Matsushita (\cite{_Matsushita:CP^n_}).
J.-M. Hwang has shown that $B$ is a complex projective
space whenether it is smooth (\cite{_Hwang:base_}).
It was conjectured that $B$ is biholomorphic
to $\C P^n$ when it is normal
(\cite{_CMSB:adv_stud_pure_}).

Since the degenerate twistor
deformation is uniquely determined
by the cohomology class of $\eta$, and $b_2(B)=1$,
all degenerate twistor deformations belong
to a 1-dimensional holomorphic family 
determined by $\Omega+t\pi^*\eta$, with $t\in \C$.
This family is obtained as a limit of twistor families,
which explains the term ``degenerate twistor
  deformation''. Like in the usual twistor case,
a generic fiber of a degenerate twistor family
is non-algebraic. However, the Moishezon fibers are
dense in the degenerate twistor family.

Much is still unknown about the degenerate twistor
families. For example, it is unknown whether
all fibers of the degenerate twistor
deformation are K\"ahler. It is unclear
when a degenerate twistor deformation
of a Lagrangian fibration admits a section.

In this paper we prove that existence
of a section is essentially a topological condition.

\hfill

\theorem\label{_Section_intro_Theorem_}
Let $\pi:\; M \arrow B$ be a Lagrangian fibration 
on a holomorphic symplectic manifold $(M, \Omega)$,
and $S:\; B \arrow M$ its smooth section.
Then there exists a closed 
form $\eta\in \Lambda^{2,0}(B)+\Lambda^{1,1}(B)$
such that $S$ is holomorphic Lagrangian with respect to the 
complex structure $I_\eta$ induced by $\Omega+
\pi^*\eta$.

\proof
\ref{sctn}.
\endproof

\hfill

Existence of
holomorphic Lagrangian fibrations is a standard
conjecture, sometimes referred to as ``hyperk\"ahler
SYZ conjecture''; see \cite{_Verbitsky:SYZ_} for
the history and a precise formulation.
It is known that holomorphic Lagrangian fibrations
exist for all known classes of hyperk\"ahler manifolds 
(see e. g. \cite[Claim 1.20]{_Kamenova_Verbitsky:fibrations_}). 
Any Lagrangian torus on a hyperk\"ahler manifold
$M$ is a fiber of a rational Lagrangian fibration;
moreover, this fibration becomes holomorphic
after replacing $M$ with another
holomorphically symplectic birational model
(\cite{_Hwang_Weiss:Beauville_,_Greb_Lehn_Rollenske_}).

Existence
of holomorphic sections would follow if we 
construct a smooth section, which seems to be
doable in many (or all) cases. 
This would bring us closer to 
understanding the holomorphic Lagrangian 
fibrations, and (hopefully) bring new 
classification results. 

There is an obvious topological obstruction to existence
of a smooth section: a fibration with multiple fibers
cannot have sections. It is easy to construct
holomorphic Lagrangian fibrations with multiple
fibers if the general fiber is not connected.
However, this problem can be rectified by
using the Stein factorization theorem.
Further on, we shall assume that
all holomorphic Lagrangian fibrations
we consider have connected fibers.

In that case (as far as we know)
there are no examples of Lagrangian fibrations 
with multiple fibers on compact, maximal
holonomy hyperk\"ahler manifolds. We conjecture that multiple
fibers in this situation don't occur.

For K3 surface non-existence of multiple fibers
is implied by Kodaira's classification of singular elliptic fibers
(it follows directly from the
canonical bundle formula, \cite[Thm 12.1]{_BHPV:04_}).
We give a direct proof of this fact in Section
\ref{_K3:Section_}. 

Multiple fibers of Lagrangian
fibrations were classified in 
\cite{_Hwang_Oguiso:multiple_} 
and \cite{_Matsushite:multiple_}.
In particular, it is shown that
the multiplicity of a general fiber
is at most 6  (\cite[Theorem 1.1]{_Hwang_Oguiso:multiple_}).

It is clear from our construction
that for any smooth section $S$ of a
Lagrangian fibration $\pi$, there exists a unique degenerate 
twistorial deformation of $\pi$
such that $S$ is holomorphic.  We 
don't know the number of degenerate 
twistorial deformations of $\pi$ admitting
a holomorphic section; it is unknown even whether this number
is finite or infinite.

%%%%%%%%%%%%%%%%%%%%%%%%%%%%%%%%%%%%%%%%%%%%%%%%%%%%%%%%%%%%%%%%%%%%%

\section{Complex structures via symplectic forms}

%%%%%%%%%%%%%%%%%%%%%%%%%%%%%%%%%%%%%%%%%%%%%%%%%%%%%%%%%%%%%%%%%%%%%

%%%%%%%%%%%%%%%%%%%%%%%%%%%%%%%%%%%%%%%%%%%%%%%%%%%%%%%%%%%%
\subsection{C-symplectic structures}
%%%%%%%%%%%%%%%%%%%%%%%%%%%%%%%%%%%%%%%%%%%%%%%%%%%%%%%%%%%%

\definition
Let $V$ be a real vector space equipped with an operator $I\in \End(V)$,
$I^2=-\Id_V$. We say that $I$ is a {\bf complex structure operator} on $V$, and
$V$ is a {\bf complex vector space}. The operator $I$ can be understood as 
imaginary unit $\1\in \C$ acting on $V$. 
Complex-linear maps of complex vector spaces
are the same as maps which commute with the complex structure.

\hfill

Suppose $W$ is a complex vector space and $\Omega \in
\Lambda^2_{\C}W^*$ a complex symplectic form. Let us
forget about the complex structure on $W$, and consider
its underlying real vector space $W_\R$ with a 2-tensor
$\Omega \in \Lambda^2_{\R}W^*\o\C$ with complex
coefficients. Knowing just this tensor, one can uniquely
reconstruct the complex structure 
operator $I_{\Omega} \colon W_{\R} \to
W_{\R}$. In the present section we determine
which forms can occur as complex symplectic
forms for some complex structure on $W_\R$.

\hfill

%%%%%%%%%%%%%%%%%%%%%%%%%%%%%%%%%%%%%%%%%%%%%%%%%%%%%%%%%%%%
\definition
Let $V$ be an $4n$-dimensional real vector space. A
2-tensor $\Omega \in \Lambda^2V^*\o\C$ is called a {\bf
  c-symplectic form}, if for any nonzero vector $v\in V$
one has $\iota_v\Omega \neq 0 \in V^*\o\C$, and
$\ker\Omega \subset V\o\C$ has rank $2n$. A pair $(V,\Omega)$ in such situation is
called a {\bf c-symplectic vector space}.

\hfill

\remark
The kernel rank $\ker \Omega=2n$ is maximal possible 
for a form which does not vanish on real vectors.
Indeed, if $\dim_\C \ker \Omega>2n$, equivalently,
$\dim_\R \ker \Omega>4n$
this space would intersect $V\subset V\otimes \C$
of real dimension $4n$, and this is impossible.

\hfill

\definition
Let $(V,\Omega)$ be a c-symplectic vector space. A complex
structure operator $I \colon V \to V$ is called an {\bf
induced} by the c-symplectic form
$\Omega$, if it makes $\Omega$ a complex linear symplectic
form in $\Lambda^2_\C V^*$, where the complex exterior
power is taken w.~r.~t. the complex structure $I$.

\hfill

\proposition\label{_inhe_complex_stru_Proposition_}
For any c-symplectic vector space $(V,\Omega)$ there exists
a unique induced complex structure $I$ on $V$, that is,
a complex structure operator $I\in \End_\R(V)$ such that
$\Omega$ is a non-degenerate form of Hodge type (2,0)

\hfill

\proof
Since $\ker\Omega \subset V\o\C$ has maximal possible
dimension, any real vector $v \in V$ can be represented as
$v = v^{1,0} + v^{0,1}$, where $v^{0,1} \in \ker\Omega$
and $v^{1,0} \in \overline{\ker\Omega}$. Moreover, this
representation is unique, since $\ker\Omega$ contains no
real vectors and hence does not intersect its complex
conjugate subspace. One can define an operator $I$ as
multiplication by $-\sqrt{-1}$ on $\ker\Omega$ and by
$\sqrt{-1}$ on $\overline{\ker\Omega}$. Since it is
self-conjugate, it is defined over reals, it is
well-defined because of the existence and uniqueness of
the above decomposition, and it obviously squares to
$-\Id_V$. Moreover, one has $\Omega(Iu,v) =
\Omega(\sqrt{-1}u^{1,0} - \sqrt{-1}u^{0,1},v) =
\sqrt{-1}\Omega(u^{1,0},v) =
\sqrt{-1}\Omega(u,v)$. Therefore $I$ is an induced complex
structure; in particular, at least one induced complex
structure exists.

On the other hand, any induced complex structure must have
$\ker\Omega$ as its $-\sqrt{-1}$-eigenspace, and this
determines a complex structure operator in a unique way.
\endproof

\hfill

We obtain that any c-symplectic vector space is of nature
prescribed above: it is the underlying real space of some
complex symplectic space. In what follows, we denote the
induced complex structure of a c-symplectic vector space
$(V,\Omega)$ by $I_\Omega$.

\hfill

\definition
Let $(V,\Omega)$ be a c-symplectic vector space. A real
subspace $U \subset V$ is called {\bf c-isotropic} if for
any $u, u' \in U$ one has $\Omega(u,u') = 0$, and {\bf
  c-Lagrangian} if it is c-isotropic and is not contained
in any proper c-isotropic supersubspace.

\hfill

\proposition (\cite[Proposition 1]{Hi})\\
\label{Hitchin}
Any c-Lagrangian subspace of a c-symplectic vector space is 
preserved by the induced complex structure.

\hfill

\proof
Let $L \subset V$ be a c-Lagrangian subspace, and $u \in
L$ be any vector. Consider the linear span of the space
$L$ and the vector $I_\Omega u$, denote it by $L_u$. What
is the restriction of the form $\Omega$ on this span? By
definition of a c-isotropic subspace, one has $\Omega(u,v)
= 0$ for any $v \in L$, so it is completely determined by
the 1-form $\left(\iota_{I_\Omega
  u}\Omega\right)|_{L_u}$. Moreover, since the form
$\Omega$ has Hodge type $(2,0)$ w.~r.~t. $I_\Omega$, one
has $0 = \sqrt{-1}\cdot 0 = \sqrt{-1}\Omega(u,v) =
\Omega(I_\Omega u,v)$, and since $\Omega$ is
skew-symmetric, one has $\Omega(I_\Omega u, I_\Omega u) =
0$. Therefore the form $\Omega$ vanishes identically on
the subspace $L_u$, i.~e. $L_u$ is c-isotropic, and since
$L$ is contained in no proper c-isotropic supersubspace,
it must be equal to $L$. This implies that $I_\Omega u \in
L$, and since $u$ is arbitrary, it means that the operator
$I_\Omega$ maps the subspace $L$ to itself.
\endproof

\hfill

In particular, quotient by a c-Lagrangian vector subspace
inherits the complex structure. Note that this complex
structure on the quotient is not in general determined by
any c-symplectic form: indeed, the quotient may have odd
complex dimension. In what follows, we shall refer to the
complex structures on quotients of c-symplectic vector
spaces by their c-Lagrangian subspaces which makes the
projection map complex linear as to {\bf inherited} ones.

\hfill

\proposition (existence of a c-symplectic
    basis).\label{Gram-Schmidt}
Let $V$ be a c-symplectic space of real dimension $4n$. 
Then $V$ possesses a basis in which its c-symplectic form
is given by a block diagonal matrix with equal blocks $Q$
on the diagonal, where $Q$ stands for a $4\times 4$-block
given by $$Q=\begin{pmatrix}0 & 0 & 1 & \sqrt{-1} \\ 0 & 0
& \sqrt{-1} & -1 \\ -1 & -\sqrt{-1} & 0 & 0 \\ -\sqrt{-1}
& 1 & 0 & 0\end{pmatrix}.
$$
\proof
The Proposition is proved by running an analogue of the
symplectic Gram--Schmidt process. Indeed, $I_\Omega$ be
the induced complex structure, and
$z_1, ..., z_{2n}\in V$ a basis in $(V^*,I)$, considered
as a complex vector space, such that
$\Omega=\sum_{i=1}^nz_{2i-1}\wedge z_{2i}$.
Then $z_1, I_\Omega(z_1), ..., z_{2n},I_\omega(z_{2i})$
is a real basis in $V$ such that $\Omega$ is written
as the above block matrix. 

For the convenience of the reader, we give a direct
argument constructing such a basis explicitly.

For the first two vectors in the basis we may choose an
arbitrary nonzero vector $u_1$ and its image under the
induced complex structure, $I_\Omega u_1$. Orthogonals of
these vectors (i.~e. the kernels of the forms
$\iota_{u_1}\Omega$ and $\iota_{I_\Omega u_1}\Omega$)
coincide by the definition of the induced complex
structure. If one picks a vector $u_2$ outside this
orthogonal, then, by outstretching it and adding to it a
multiple of the vector $I_\Omega u_2$ in case of
necessity, we can make $u_2$ such
that $$\Omega(u_1,u_2)=1.$$ One can deduce from complex
linearity of $\Omega$ w.~r.~t. $I_\Omega$ and the above
relation that the restriction of the form $\Omega$ onto
the four-dimensional subspace $U\subset V$ spanned by
$\{u_1,I_\Omega u_1,u_2,I_\Omega u_2\}$ is given by the
matrix $Q$ in these coordinates.

Now we can proceed, exercising the same procedure in the
orthogonal to $U$ (i.~e. the set $U^\perp$ of vectors $w
\in V$ s.~t. $\left(\iota_w\Omega\right)|_U = 0$), because
$U\cap U^\perp=\{0\}$. 
\endproof

\hfill

\claim\label{_Omega^n+1=0_equiv_Claim_}
Let $V$ be a real vector space of real dimension $4n$, and
$\Omega \in \Lambda^2V^*\o\C$ be a complex-valued
skew-symmetric 2-form. The following are equivalent:
\begin{enumerate}
	\item $\dim_\C\ker\Omega\restrict{V\otimes \C}=2n$,
	\item $(\Omega\wedge\overline{\Omega})^n$ is 
nonzero and $\Omega^{n+1}=0$.
\end{enumerate}

\proof
Suppose that $\dim\ker\Omega=2n$. Then for any
decomposable polyvector $a\in\Lambda^{2n+2}V\otimes\C$ the
corresponding subspace in $V$ needs to intersect
$\ker\Omega$, so $\iota_a\Omega^{n+1}$ vanishes. Hence
$\Omega^{n+1}=0$. An easy direct calculation shows that,
in notation of the \ref{Gram-Schmidt},
$(Q\wedge\overline{Q})(u_1,v_1,u_2,v_2)=4$, so the top
power of $\Omega\wedge\overline{\Omega}$ cannot be zero.

Suppose that $\Omega^{n+1}=0$. If the top power of a
skew-symmetric form on some vector space is zero, then
this form has a nontrivial kernel. That's why the
form~$\Omega$ has nontrivial kernel when restricted to any
$(2n+2)$-plane inside $V\otimes\C$ (and, moreover, any
$(2n+2k)$-plane for any $k>0$). Restriction of a
symplectic form onto a maximal subspace tranversal to its
kernel is non-degenerate, so $\dim\ker\Omega\geqslant
2n$. If $(\Omega\wedge\overline{\Omega})^n \neq 0$, it
cannot be greater, because in this case $\ker\Omega$ needs
to intersect $\ker\overline{\Omega}$, thus giving a real
vector in $\ker\Omega$, substitution of which would vanish
$(\Omega\wedge\overline{\Omega})^n$.
\endproof

% Is this proposition even needed? -- RD
%
%\begin{pr}
%\label{_antilinear_Claim_}
%Let $(A,I)$, $(B,I)$  be complex vector spaces,
%and $\lambda:\; A \to B$ a real map. Then $\mu:=[\lambda, I]$ is complex antilinear,
%that is, satisfies $[I, \mu]=-I$.
%\end{pr}
%\begin{proof}
%Clear.
%\end{proof}

%%%%%%%%%%%%%%%%%%%%%%%%%%%%%%%%%%%%%%%%%%%%%%%%%%%%%%%%%%%%
\subsection{Sections of Lagrangian projections of 
c-symplectic vector spaces}
%%%%%%%%%%%%%%%%%%%%%%%%%%%%%%%%%%%%%%%%%%%%%%%%%%%%%%%%%%%%

In this subsection, we will
deal with the following situation.
Let $(V,\Omega)$ be a $4n$-dimensional 
real vector space space equipped with
 a c-symplectic structure, and
$L\subset V$ a c-Lagrangian subspace.
The complex linear surjective map $V\arrow V/L$ is
called {\bf a c-Lagrangian projection}.
This is a linearization of the Lagrangian
fibrations which are common in
holomorphic symplectic geometry
(Subsection \ref{_Lagra_fibra_Subsection_}).

\hfill

\proposition
\label{calculation}
Let $V$ be a real vector space, $\Omega$ a c-symplectic
form on it, $L$ a c-Lagrangian subspace and $\sigma \colon
V/L \to V$ a real section (not necessarily complex linear
w.~r.~t. the induced complex structures). Define the form
$\Omega_\sigma\in\Lambda^2(V/L)^*\otimes\C$ by the rule
$\Omega_\sigma(u_1,u_2) =
\Omega(\sigma(u_1),\sigma(u_2))$.
In other words, $\Omega_\sigma$ is a 
restriction of $\Omega$ onto the subspace
$\sigma(V/L)$ after the identification
$\pi|_{\sigma(V/L)}\colon\sigma(V/L)\to V/L$). Then 
$\Omega_\sigma$ has Hodge type $(2,0)+(1,1)$
w.~r.~t. the inherited complex structure on the quotient
$V/L$.

\hfill

\proof
Suppose that $\sigma = \sigma_0$ is complex linear. Then
the form $\Omega_\sigma$ is of type $(2,0)$, since the
type is preserved by complex linear maps.

Now let $\sigma = \sigma_0 + \tau$, where $\sigma_0$ is a
complex linear section and $\tau \colon V/L \to L$ a
real perturbation. For any vectors $u,v\in (V/L)$,
\begin{multline} \label{_restrict_to_section_2,0+1,1_Equation_}
\Omega_\sigma(u,v) = \Omega(\sigma(u),\sigma(v)) = \\=
\Omega(\sigma_0(u),\sigma_0(v)) +
\Omega(\sigma_0(u),\tau(v)) +\\+
\Omega(\tau(u),\sigma_0(v))  + \Omega(\tau(u),\tau(v)).
\end{multline}
Since one has $\tau(u),\tau(v) \in L$ and $L$ is a
c-Lagrangian subspace w.~r.~t. $\Omega$, the last term
vanishes. The first term is a $(2,0)$ form, since
$\sigma_0$ is a complex linear section. The term
$\Omega(\sigma_0(u),\tau(v)) +
\Omega(\tau(u),\sigma_0(v))$ vanishes if $u,v$ are both of
type $(0,1)$, since $\sigma_0(u)$ and $\sigma_0(v)$ are
vectors of type $(0,1)$ in this case, and therefore
annihilate the form $\Omega$. Therefore, all terms of 
\eqref{_restrict_to_section_2,0+1,1_Equation_} 
are of type $(2,0)+(1,1)$.
\endproof

\hfill

Fiberwise application of the above construction allows one
to obtain an almost complex structure on a
$4n$-dimensional manifold $X$ from a non-degenerate
complex-valued 2-form $\Omega$ on it such that
$(\Omega\wedge\overline{\Omega})^{n}$ is nowhere zero and
$\Omega^{n+1}=0$; by \ref{_Omega^n+1=0_equiv_Claim_}, this condition
gives a c-symplectic form on the tangent bundle. 

\hfill

\definition
An {\bf almost c-symplectic form} on a manifold $X$ of
real dimension $4n$ is a form
$\Omega\in\Gamma(\Lambda^2T^*X\otimes\C)$ such that
the~top degree form $(\Omega\wedge\overline{\Omega})^n$ is
nowhere zero and $\Omega^{n+1}=0$. A {\bf c-symplectic
  form} is a closed almost c-symplectic form.

\hfill

We shall denote the induced almost complex structure
obtained from the form $\Omega$ as the induced complex
structure on each fiber by $I_\Omega$.

\hfill

\proposition
[\cite{_Verbitsky:Degenerate_} , Theorem 3.5]\label{clos}\\
Let $\Omega$ be a c-symplectic form on a manifold,
and $I_\Omega\in \End(TM)$ the corresponding almost complex structure.
Then $I_\Omega$ is integrable. 

\hfill

\proof
From Cartan's formula it follows immediately that
for any closed $k$-form $\Phi$ and two vectors
$X, Y$ such that the contractions of $\Phi$ with $X, Y$
vanish, $\Phi\cntrct X = \Phi\cntrct Y=0$, one also has
$\Phi\cntrct ([X,Y])=0$. However, the space
$T^{0,1}_{I_\Omega}M\subset TM\otimes \C$
is defined as
\[
T^{0,1}_{I_\Omega}M=\{ X\in TM\otimes
\C\ \ |\ \ \Omega\cntrct X=0\},
\]
which gives $[T^{0,1}_{I_\Omega}M, T^{0,1}_{I_\Omega}M]\subset T^{0,1}_{I_\Omega}M$.
By Newlander-Nirenberg, this condition is equivalent
to integrability of $I_\Omega$.
\endproof

Note that converse is not generally true: if $\Omega$ is a
c-symplectic form and $f$ is a non-vanishing function, then
$f\Omega$ is not closed unless $f$ is constant; however,
$I_{f\Omega}=I_\Omega$.

%%%%%%%%%%%%%%%%%%%%%%%%%%%%%%%%%%%%%%%%%%%%%%

\section{Degenerate twistorial deformation}

%%%%%%%%%%%%%%%%%%%%%%%%%%%%%%%%%%%%%%%%%%%%%%

Though this section is self-contained, its results
generalise and simplify the results of
\cite{_Verbitsky:Degenerate_}.

\hfill

\proposition
\label{degenerate_linear}
Let $\Omega\in\Lambda^2V^*\otimes\C$ be a
c-symplectic form on a real vector space $V$, $L$ a
c-Lagrangian subspace in $V$ and $\pi\colon V\to V/L$ the
projection. Then for any  complex valued form
$\gamma\in\Lambda^2(V/L)^*\otimes\C$ of Hodge type
$(2,0)+(1,1)$ w.~r.~t. the inherited complex structure on
$V/L$, the form $\Omega_\gamma=\Omega+\pi^*\gamma$ is
c-symplectic.

\hfill

\proof
We shall proceed in two steps: first, we prove that
$\dim\ker\Omega_\gamma$ is at least half of $\dim V$, and
then we prove that $\dim\ker\Omega_\gamma$ does not exceed
half of $\dim V$. 

\hfill

{\bf Assertion 1: $\dim\ker\Omega_\gamma\geqslant \frac12\dim
  V$.} 
Clearly, this would follow if we prove that $\Omega_\gamma^{n+1}=0$,
where $4n=\dim_\R V$.
Let $V_\C:= V\otimes_\R \C$ be the complexification
of $V$, $L\oplus K = V$ a direct sum decomposition,
with $K,L\subset V$ complex vector spaces,
and $L_\C, K_\C$ complexifications of these spaces. Denote by
$L_\C= L_\C^{1,0}\oplus  L_\C^{0,1}$ and  
$K_\C= K_\C^{1,0}\oplus  K_\C^{0,1}$  their Hodge decompositions.
Consider a linear automorphism 
$T_\lambda:\;V_\C \arrow V_\C$ acting as multiplication by a
scalar $\lambda\in \C$ on 
$K_\C^{1,0}$ and as identity on $L_\C$ and $K_\C^{0,1}$.
Since $\Omega$ is a pairing between $K_\C^{1,0}$ and
$L_\C^{1,0}$, one has $T_\lambda(\Omega)=\lambda\Omega$.
Since $\pi^*\gamma$ vanishes on $L$, one has 
$T_\lambda(\pi^*\gamma) = \lambda^2(\pi^*\gamma)^{2,0} + \lambda(\pi^*\gamma)^{1,1}$.
Using this, we write the weight decomposition for the action of $T_\lambda$ on
$\Omega_\gamma^{n+1}$ as follows:
\[ 
  T_\lambda(\Omega_\gamma^{n+1})=
  \sum_{i=0}^n \sum_{j=0}^{n-i}\lambda^{n+1+j}\Omega^i\wedge (\gamma^{2,0})^j \wedge 
(\gamma^{1,1})^{n-i-j+1}
\]
However, $T_\lambda$ cannot act on $2n+2$-forms with weight $\geq n+1$
because $\dim_\C L_\C^{1,0}=n$. Therefore, $\Omega_\gamma^{n+1}=0$.

\hfill

{\bf Assertion 2: $\dim\ker\Omega_\gamma\leqslant \frac12\dim
  V$, or, equivalently, $\rk \Omega_\gamma\geq
  \rk\Omega$.} 
Clearly, the point $\Omega \in \Lambda^2(V_\C)$ 
has a neighbourhood $U$ such that all 2-forms $y\in U$
have rank $\geq \rk \Omega$. Assertion 2 is trivial
when $\Omega_\gamma\in U$. Consider, as in the previous step,
the decomposition 
\[ V_\C=L_\C^{1,0}\oplus  L_\C^{0,1}\oplus K_\C^{1,0}\oplus  K_\C^{0,1},
\]
and let $R_\lambda\in \Aut(V_\C)$ act on $K_\C^{0,1}$ and
$L_\C^{0,1}$ as identity,
on $L_\C^{1,0}$ as multiplication by $\lambda^{-1}$ and on $K_\C^{1,0}$
as multiplication by $\lambda$. Since $\Omega$ is a pairing
between $L_\C^{1,0}$ and $K_\C^{1,0}$, one has $R_\lambda(\Omega)=\Omega$.
Since $\gamma$ vanishes on $L$, one has 
$R_\lambda(\gamma) = \lambda^2\gamma^{2,0} + \lambda\gamma^{1,1}$.
Then $\lim_{\lambda\mapsto 0}R_\lambda(\Omega_\gamma)=\Omega$,
hence for $\lambda$ sufficiently small, the form 
$R_\lambda(\Omega_\gamma)$ belongs to $U$, and 
satisfies $\rk \Omega_\gamma \geq \rk \Omega$.
We proved \ref{degenerate_linear}. \endproof

\hfill

\proposition
\label{preservance}
Let $\Omega\in\Lambda^2V^*\otimes\C$ be a
c-symplectic form on a real vector space $V$, $L$ a
c-Lagrangian subspace in $V$ and $\pi\colon V\to V/L$ the
projection. Consider a (2,0)+(1,1)-form $\gamma$ on
$V/L$, and let $\Omega_t:=\Omega+t\pi^*\gamma$, for some
$t\in \C$. Then the subspace $L$ is c-Lagrangian
w.~r.~t. all the forms $\Omega_t$, and restrictions
$I_{\Omega_t}\restrict L$ of induced complex structures
$I_{\Omega_t}$ coincide. Moreover, the complex
structures on the quotient $V/L$ inherited from $I_{\Omega_t}$
also coincide, for all $t\in \C$.

\hfill

\proof
The first claim is obvious from the construction of $\Omega_t$.
The prove the second claim, take
$v\in(V/L)^{0,1}_{I_\Omega}$, and 
let $u\in \pi^{-1}(v)$ be a vector in its preimage.
Then $\Omega_t\cntrct u = t\pi^*(\gamma\cntrct v)$.
Therefore, 
\begin{equation} \label{_varies_0,1_Equation_}
u\in V^{0,1}_{I_{\Omega_t}}\ \  \Leftrightarrow\ \ \forall
z\in V, \ \ 
\Omega(u, z)= - t\pi^*(\gamma\cntrct v)(z)) 
\end{equation} The 1-form 
$z\arrow t\pi^*(\gamma\cntrct v)(z)$
vanishes on $L$ and has type (1,0)
on $(V, I_\Omega)$, because 
$\gamma$ is of type (2,0)+(1,1).

Since $\Omega$ is a non-degenerate
pairing between $L^{0,1}_{I_\Omega}$ and $(V/L)^{0,1}_{I_\Omega}$, 
for any (1,0)-form $\xi$ on $V/L$
there exists a vector $x\in L$
such that $\xi= \Omega\cntrct x$. 

This gives a vector $\zeta_t\in L^{0,1}_{I_\Omega}$ such that
such that \[ \Omega(u+\zeta_t, z)= - t\pi^*(\gamma\cntrct v)(z))\]
for all $z\in V$. By \eqref{_varies_0,1_Equation_}, 
$u+\zeta_t\in V^{0,1}_{I_{\Omega_t}}$
is a vector which projects to $v$.

This implies that the space of (0,1)-vectors
in $(V/L)_{I_{\Omega_t}}$ is independent from $t$.
\endproof

\hfill

\theorem
\label{degenerate}
Suppose that $X\arr{\pi}B$ is a Lagrangian fibration on a
holomorphically symplectic manifold $(X,\Omega)$, and
$\eta\in\Omega_{\mathrm{cl}}(B)$ a closed
$(2,0)+(1,1)$-form on the base. Then the forms
$\Omega_t=\Omega+t\pi^*\eta$ on $X$ are c-symplectic, and
this deformation (called {\bfseries\itshape degenerate
  twistorial deformation}) preserves the Lagrangian
fibration and the base.

\hfill

\proof
From \ref{degenerate_linear} it follows
that $\Omega_t$ is c-symplectic. 
The fibers stay Lagrangian, and the complex
structure thereof, as well as that of the base, remains
unchanged by \ref{preservance}.
\endproof

\hfill

\remark This result was proven for compact manifolds in 
\cite[Theorem 1.10]{_Verbitsky:Degenerate_}.

\hfill

For a hyperk\"ahler manifold $X$, the degenerate
twistorial deformation produces an entire curve in the
period space of $X$. In terms of the oriented 2-plane
Grassmannian (see e.~g. Section 3 in \cite{De}), the plane
corresponding to $(X,\Omega_{x+\sqrt{-1}y})$ is spanned by
$(\Omega+\overline{\Omega}) + 2x\eta$ and
$\sqrt{-1}(\Omega-\overline{\Omega}) - 2y\eta$. Thus one
can define degenerate twistorial curves in an abstract
situation: namely, for a 2-plane $\tau\in\Gr_{++}(V,q)$
and a vector $e\in \tau^\perp \subset V$ with $q(e,e)=0$,
the subvariety
$\Deg_\tau(e)=\Gr_{++}(\mathrm{span}(\tau,e),q|_{\mathrm{span}(\tau,e)})\subseteq\Gr_{++}(V,q)$
is an entire curve, and in the case when $V=H^2(X,\R)$,
$q$ is the Bogomolov--Beauville--Fujiki form and
$e=\pi^*[\omega]$ is the inverse image of the K\"ahler
class on the base, this curve is exactly the base of the
degenerate twistorial deformation.

\hfill

\theorem
\label{sctn}
Let $X \arr{\pi} B$ be a holomorphic Lagrangian fibration on a
holomorphically symplectic manifold $(X,\Omega)$, and
$\sigma \colon B \to X$ a smooth section of $\pi$. Then there exists a
degenerate twistorial deformation $(X,\Omega')$ of $(X, \Omega)$
s.~t. the fibers stay Lagrangian, the complex
structure on the fibers and the base stays the same, 
and $\sigma$ is a holomorphic map.

\hfill

\proof
Consider the form $\eta = \sigma^*\Omega \in
\Omega^2(B)$. By Proposition \ref{calculation}, it is of
type $(2,0)+(1,1)$. Then by Proposition \ref{degenerate}
the forms $\eta$ gives rise to a deformation with desired
properties. For $t=-1$ one has $\Omega_t|_{\sigma(B)} =
\left(\Omega - \pi^*\sigma^*\Omega\right)|_{\sigma(B)} =
\Omega|_{\sigma(B)} - \Omega|_{\sigma(B)} = 0$. By
Hitchin's lemma (\ref{Hitchin}), this means that the
submanifold $\sigma(B)$ is a complex submanifold. 
Since the projection is a holomorphic map, the section
$\sigma$ is also holomorphic.
\endproof

%%%%%%%%%%%%%%%%%%%%%%%%%%%%%%%%%%%%%%%%%%%%%%%

\section{Holomorphic Lagrangian fibrations on a K3
  surface}
\label{_K3:Section_}

%%%%%%%%%%%%%%%%%%%%%%%%%%%%%%%%%%%%%%%%%%%%%%%

In this section we apply
\ref{sctn} to obtain 
sections of holomorphic Lagrangian fibrations on
K3 surfaces. 

It is not hard to see that the elliptic
fibrations on a K3 surface 
$X$ are in bijective correspondence with primitive numerically
effective classes $e\in\NS(X)$ with $(e,e)=0$.
Indeed, its linear system has no basepoints and
establishes a map to $\C P^1$ with generic fiber elliptic
curve. Cross-sections of the fibration determined by $e$
correspond to classes $s\in\NS(X)$ with $(s,e)=1$. Our
first goal is to find an effective class $s\in H^2(X,\Z)$
with $(s,e)=1$ and a deformation in which $s$ would have
type $(1,1)$.

\hfill

\claim
\label{unimod}
For any isotropic primitive vector $e$ in an even
unimodular lattice $\Lambda$ there exists a vector
$a\in\Lambda$ such that one has $(a,e)=1$ and $(a,a)=-2$.

\hfill

\proof
Since the lattice $\Lambda$ is unimodular, one can pick
some vector $b\in\Lambda$ with $(b,e)=1$. The number
$(b,b)$ is even since $\Lambda$ is even, so the vector
$a=b-\left(1-(b,b)/2\right)e$ is integral. One has $(a,e)
= (b,e) - (1+{(b,b)}/2)(e,e) = 1$ and $(a,a) =
(b-(1+{(b,b)}/2)e,b-(1+{(b,b)}/2)e) = (b,b) -
2(1+{(b,b)}/2)(b,e) = (b,b) - 2 - (b,b) = -2$.
\endproof

\hfill

The following claim is trivial.

\hfill

\claim\label{_FS_fiber_Claim_}
Suppose that $X$ is a K3 surface with an elliptic
fibration determined by an isotropic class $e$. Then the
pullback of the Fubini--Study form from the base
represents $e$.

\hfill

\proof
The Fubini--Study form on the base represents a Poincar\'e
dual to the class of a point. Then its inverse image
represents the Poincar\'e dual to the class of the inverse
image of the point, i.~e. of the fiber.
\endproof

\hfill

\lemma
\label{construction_k3}
Suppose that $X$ is a K3 surface with elliptic fibration
$X\arr{\pi}\C P^1$ determined by an isotropic class $e$,
$\omega_{\FS}$ the Fubini--Study form on its base, and
$s\in H^2(X,\Z)$ is such that one has $(s,e)\neq 0$. Then
there exists a unique degenerate twistorial
deformation $X'$ of $(X,\pi)$ in such way
that $s$ belongs to $H^{1,1}(X')$.

\hfill

\proof
The forms $\Omega-t\pi^*\omega_{\FS}$, $t\in\C$, are
c-symplectic by Proposition \ref{degenerate}.
Denote by $I_t$ the complex structure induced by 
$ \Omega-t\pi^*\omega_{\FS}$ (\ref{_inhe_complex_stru_Proposition_}).
The homology class of $s$ is of type $(1,1)$ if and only if
$(s,[\Omega-t\pi^*\omega_{\FS}])=0$. Using
\ref{_FS_fiber_Claim_},
 one can rewrite this equation as
$(s,[\Omega])=t(s,e)$. It has a unique solution $t$ whenever
$(s,e)\neq 0$.
\endproof

\hfill

\proposition
\label{k3case}
Let $X$ be a K3 surface with elliptic fibration
$X\arr{\pi}\C P^1$. Then
there exists a degenerate twistorial
deformation $X'$ of $(X,\pi)$ admitting a 
holomorphic section.

\hfill

\proof
Since the lattice $H^2(X,\Z)$ for a K3 surface $X$ is even
and unimodular, by \ref{unimod}, one can find a vector
$s\in H^2(X,\Z)$ such that $(s,e)=1$ and $(s,s)=-2$. By
\ref{construction_k3}, we can deform the complex structure
on $X$ in such way that $s$ would have type $(1,1)$. By
Lefschetz theorem on $(1,1)$-classes, there exists a line
bundle $L\to X$ with $c_1(L)=s$. Riemann--Roch theorem for
this bundle $L$ reads $\chi(L) = \chi(\sO_X) +
\frac{L\cdot(L\otimes K_X^*)}2$. By Serre's duality one
has $h^2(L) = h^0(L^*\otimes K_X)$, and, as soon as $K_X$
is trivial and $\chi(\sO_X)=2$, we can rewrite it as
$h^0(L) - h^1(L) + h^0(L^*) = 2 + \frac{-2}2 = 1$. Since
$L$ is nontrivial, either $h^0(L)$ or $h^0(L^*)$ does not
vanish, and either $L$ or $L^*$ is effective. If $L^*$ is,
then $-s$ is represented by a curve, and $0<(-s,e) =
-(s,e)=-1$, hence $L$ is effective and $s$ is represented
by a curve. Suppose that $s=\sum_is_i$, where $s_i$ are
the classes of irreducible curves. One has
$1=(s,e)=(\sum_is_i,e)=\sum_i(s_i,e)$. All the numbers
$(s_i,e)$ are positive integers with sum $1$, so exactly
one of them, say $(s_0,e)$, equals $1$. The curve
represented by the class $s_0$ intersects each fiber at
one point, i.~e. is a section of the fibration.
\endproof

\hfill

In particular, this implies that 
elliptic fibrations on K3 surfaces cannot have multiple fibers.

\paragraph*{Acknowledgements.}
We are grateful to Dmitri Kaledin for insightful discussions.

\noindent {\sc Fedor A. Bogomolov\\
Department of Mathematics\\
Courant Institute, NYU \\
251 Mercer Street \\
New York, NY 10012, USA,} \\
\tt bogomolov@cims.nyu.edu, also: \\
{\sc National Research University, Higher School of Economics, Moscow, Russia,}
\\

\noindent {\sc Rodion N. D\'eev\\
Independent University of Moscow\\
119002, Bolshoy Vlasyevskiy Pereulok 11\\ 
Moscow, Russian Federation} \\
\tt deevrod@mccme.ru, also:\\
{\sc Department of Mathematics\\
Courant Institute, NYU \\
251 Mercer Street \\
New York, NY 10012, USA,} \\
\tt rodion@cims.nyu.edu
\\

\noindent {\sc Misha Verbitsky\\
{\sc Instituto Nacional de Matem\'atica Pura e
              Aplicada (IMPA) \\ Estrada Dona Castorina, 110\\
Jardim Bot\^anico, CEP 22460-320\\
Rio de Janeiro, RJ - Brasil }\\
also:\\
{\sc Laboratory of Algebraic Geometry,\\
National Research University HSE,\\
Department of Mathematics, 7 Vavilova Str. Moscow, Russia,}\\
\tt  verbit@mccme.ru}.

\end{document}